\input amstex
\documentstyle{amsppt}
\pagewidth{6.5in}
\pageheight{9in}
\topmatter
\title A Direct Sum decomposition for dual spaces.\endtitle
\author Javier H. Guachalla\endauthor
\address Universidad Mayor de San Andr\'{e}s.\\
La Paz Bolivia \endaddress
\email jguachal@entelnet.bo\endemail
\date October, 2005 \footnote{jguachal\@entelnet.bo} \enddate
\abstract
We developpe a direct sum decomposition for n-dual spaces.
\endabstract
\endtopmatter

\document
\head Introduction \endhead

Let $X$ be a Banach space. We will denote $\pi _{0}:X\rightarrow X^{**}$ and 
$\pi _{1}:X^{*}\rightarrow X^{***}$ the canonical maps of $X$ and its dual
space $X^{*}$ into their respective bidual spaces. According to \cite{G} we
have $X^{***}=Ran(\pi _{1})\oplus Ker(\pi _{0}^{*})$, where $\pi _{0}^{*}$
denotes the adjoint map of $\pi _{0}$.

\remark{1. Remark} In \cite{G} we also saw as an application that $X$ is reflexive if and only if its dual space $X^{*}$ is reflexive. And since the dual space of the
latter is the bidual space of $X$, in thiis case we conclude that $X^{**}$ is reflexive also. And so on, we may induce that $X$ is reflexive if and only if for all $n$ the $n$th-dual space is reflexive.
\endremark

As the reader can see the notation becomes cumbersome when we consider the $n$th-dual space.
$$X^{*\cdots *}\qquad \text{with }*\cdots *\text{ }n\text{-times}$$
And since the rest of the paper deals with these, we will use the following
notation. Let $\Bbb{N}$ denote the natural numbers with $0$ included.

\definition{2. Definition} Let $X$ be a Banach space. For $n\in \Bbb{N},$ $n\geq 1$   denote $X^{*n}$ the $n$th$-$dual space of $X$. Then we will denote its canonical map into its bidual space 
$$\pi _{n}:X^{*n}\rightarrow X^{*(n+2)}$$
with its adjoint map
$$\pi _{n}^{*}:X^{*(n+3)}\rightarrow X^{*(n+1)}$$
\enddefinition
Note that these dual spaces are all Banach spaces, and that we have the maps
$$\pi _{n}:X^{*n}\rightarrow X^{*(n+2)}\qquad \pi
_{n-1}^{*}:X^{*(n+2)}\rightarrow X^{*n}$$
and since the cited result is valid for them we state the result for sake of
completeness as

\proclaim{3. Proposition}
For $n\in \Bbb{N},$ $n\geq 1.$
$$X^{*(n+2)}=Ran(\pi _{n})\oplus Ker(\pi _{n-1}^{*})$$
\endproclaim

\head A direct sum decomposition for dual spaces\endhead

\definition{4. Definition}
Two Banach spaces $Y,$ $M$ will be said isomorphic if there exist a bijective linear operator, which, as well as its inverse are bounded, what we will denote $Y\cong M$.
\enddefinition
\remark{5. Remark}
Let $X$, $Y$, $Z$ be Banach spaces. If there exist closed subspaces $M$ and $N$ of $X$, such that $X=M\oplus N$, and $Y\cong M$ and $Z\cong N$. Where $X=M\oplus N$ means that $X=M+N$ and $M\cap N=\{0\}$. Then $X$ is isomorphic to the direct sum of $Y$ and $Z$, denoted $X\cong Y\oplus Z$, 
\endremark
\proclaim{6. Corollary}
For all $n\geq 1,$ it is valid 
$$X^{*(n+2)}\cong X^{*n}\oplus Ker(\pi _{n-1}^{*})\qquad \tag * $$
\endproclaim
\demo{Proof}
This is because $X^{*n}\cong Ran(\pi _{n})$, isomorphism which is isometric,
too.
\enddemo

Therefore, when $X$ is reflexive, the following holds

\proclaim{7. Corollary}
If $X$ is reflexive. Then, for $n\ge 1$
$$
\gather 
X \cong X^{**}\cong X^{*4}\cong X^{*6}\cong \cdots \cong X^{*(2n)}\cong \cdots \\
X^{*} \cong X^{*3}\cong X^{*5}\cong X^{*7}\cong \cdots \cong X^{*(2n+1)}\cong \cdots
\endgather
$$
\endproclaim

\proclaim{8. Corollary}
If for an $i\in \Bbb{N} $, $i\geq 0$ $Ker(\pi _{i}^{*})=0$, then $\forall
n\in \Bbb{N} $ $Ker(\pi _{n}^{*})=0$
\endproclaim

For the general case, we obtain

\proclaim{9. Theorem}
For $n\in \Bbb{N} $ $n\geq 1$ we have
$$
\align
X^{*(2n+1)} &\cong X^{*}\oplus Ker(\pi _{0}^{*})\oplus Ker(\pi _{2}^{*})\oplus
\cdots \oplus Ker(\pi _{2n-2}^{*}) \\
X^{*(2n+2)} &\cong X^{**}\oplus Ker(\pi _{1}^{*})\oplus Ker(\pi
_{3}^{*})\oplus \cdots \oplus Ker(\pi _{2n-1}^{*})
\endalign
$$
\endproclaim
\demo{Proof}
Writing a few first cases of $(*)$, we obtain

$$
\align
X^{*3} &\cong X^{*}\oplus Ker(\pi _{0}^{*}) \\
X^{*4} &\cong X^{**}\oplus Ker(\pi _{1}^{*}) \\
X^{*5} &\cong X^{*3}\oplus Ker(\pi _{2}^{*}) \\
X^{*6} &\cong X^{*4}\oplus Ker(\pi _{3}^{*}) \\
X^{*7} &\cong X^{*5}\oplus Ker(\pi _{4}^{*})
\endalign
$$

Which by isomorphism leads to 
$$
\gather
X^{*6} \cong X^{**}\oplus Ker(\pi _{1}^{*})\oplus Ker(\pi _{3}^{*}) \\
X^{*7} \cong X^{*3}\oplus Ker(\pi _{2}^{*})\oplus Ker(\pi _{4}^{*})\cong
X^{*}\oplus Ker(\pi _{0}^{*})\oplus Ker(\pi _{2}^{*})\oplus Ker(\pi _{4}^{*})
\endgather
$$

Therefore for a proof by induction we admit these true for $n=k$, that is

$$
X^{*(2k+1)}\cong X^{*}\oplus Ker(\pi _{0}^{*})\oplus Ker(\pi _{2}^{*})\oplus
\cdots \oplus Ker(\pi _{2k-2}^{*})
$$
and
$$
X^{*(2k+2)}\cong X^{**}\oplus Ker(\pi _{1}^{*})\oplus Ker(\pi
_{3}^{*})\oplus \cdots \oplus Ker(\pi _{2k-1}^{*})
$$

And prove for $n=k+1.$

Since $X^{*(m+2)}\cong X^{*m}\oplus Ker(\pi _{m-1}^{*}),$ for $m\geq 1$.
Let us take the value $m=2k+1$, then

$$
X^{*(2(k+1)+1)}\cong X^{*(2k+1)}\oplus Ker(\pi _{2k}^{*})
$$

And by hypothesis of induction

$$
X^{*(2(k+1)+1)}\cong (X^{*}\oplus Ker(\pi _{0}^{*})\oplus Ker(\pi
_{2}^{*})\oplus \cdots \oplus Ker(\pi _{2k-2}^{*}))\ \oplus Ker(\pi _{2k}^{*})
$$
which is the expresion for n=k+1, in the odd case.
Similarly, for $m=2k+4$ we have

$$
X^{*(2(k+1)+2)}\cong X^{*(2k+2)}\oplus Ker(\pi _{2k+1}^{*})
$$

Then

$$
X^{*(2(k+1)+2)}\cong (X^{*}\oplus Ker(\pi _{1}^{*})\oplus Ker(\pi
_{3}^{*})\oplus \cdots \oplus Ker(\pi _{2k-1}^{*}))\ \oplus Ker(\pi
_{2k+1}^{*})
$$
\enddemo


\Refs
\ref \key G
\by  Guachalla J
\paper A Characterization of Reflexivity for Dual Banach Spaces
\jour in ArXiv.org/math/FA0509683
\yr 2005
\endref
\endRefs
\end{document}